\documentclass[9pt,final,twocolumn]{IEEEtran}

\ifCLASSINFOpdf
\else
\fi
\usepackage{cite}
\usepackage{amssymb}
\usepackage{amsmath}
\usepackage{graphicx}
\newtheorem{assumption}{Assumption}
\newtheorem{lemma}{Lemma}
\newtheorem{theorem}{Theorem}

\newtheorem{remark}{Remark}

\newtheorem{corollary}{Corollary}

\usepackage{epstopdf}
\usepackage{color}
\usepackage{cases}

\usepackage{subfigure}

\usepackage{bm}

\usepackage{algorithm}
\usepackage{algorithmicx}

\begin{document}
%
\title{Decentralized Online Learning for Noncooperative Games in Dynamic Environments}
%
%
%

\author{
Min Meng, Xiuxian Li, {\em Member, IEEE}, Yiguang Hong, {\em Fellow, IEEE}, Jie Chen, {\em Fellow, IEEE}, and Long Wang
\thanks{This work was partially supported by the National Natural Science Foundation of China under Grant 62003243 and 61733018, Shanghai Municipal Commission of Science and Technology No. 19511132101, and by Shanghai Municipal Science and Technology Major Project under grant 2021SHZDZX0100. (\emph{Corresponding author: Xiuxian Li.})}
\thanks{M. Meng, X. Li, Y. Hong, and J. Chen are with the Department of Control Science and Engineering, College of Electronics and Information Engineering, and Shanghai Research Institute for Intelligent Autonomous Systems, Tongji University, Shanghai, China (email:  mengmin@tongji.edu.cn); xli@tongji.edu.cn; yghong@tongji.edu.cn; chenjie206@tongji.edu.cn).

M. Meng and X. Li are also with the Institute for Advanced Study, Tongji University, Shanghai, China.

L. Wang is with the Center for Systems and Control, College of Engineering, Peking University, Beijing, China (e-mail: longwang@pku.edu.cn).
}
}

\maketitle

\begin{abstract}
Decentralized online learning for seeking generalized Nash equilibrium (GNE) of noncooperative games in dynamic environments is studied in this paper. Each player aims at selfishly minimizing its own time-varying cost function subject to time-varying coupled constraints and local feasible set constraints. Only local cost functions and local constraints are available to individual players, who can receive their neighbors' information through a fixed and connected graph. In addition, players have no prior knowledge of cost functions and local constraint functions in the future time. In this setting, a novel distributed online learning algorithm for seeking GNE of the studied game is devised based on mirror descent and a primal-dual strategy. It is shown that the presented algorithm can achieve sublinearly bounded dynamic regrets and constraint violation by appropriately choosing decreasing stepsizes. Finally, the obtained theoretical result is corroborated by a numerical simulation.
\end{abstract}


%
\IEEEpeerreviewmaketitle

\section{Introduction}
Game theory has received growing attention recently owing to its wide applications in social networks \cite{ghaderi2014opinion}, sensor networks \cite{stankovic2012distributed}, smart grid \cite{saad2012game}, and so on. In noncooperative games, the concept of \emph{Nash equilibrium} (NE) plays a pivotal role by providing a rigorous mathematical characterization of the stable and desirable states, from which players have no incentive to deviate \cite{basar1999dynamic}.

A challenge is to design distributed algorithms for seeking NE in noncooperative games based on limited information available to each player. Generally, it is assumed that each player can access all the players' decisions. For example, a coordinator exists to broadcast the data to the players \cite{facchinei2010generalized,yu2017distributed,shamma2005dynamic}, that is, bidirectional communication with all the agents is required, which results in high communication loads and is impractical for many applications. Therefore, distributed algorithms on computing NEs in noncooperative games without full action information, that is, in a partial decision information setting, have been getting more and more attention in recent years. To deal with such kind of scenarios, numerous results on the NE seeking problems have sprung up both in continuous-time \cite{de2019distributed,gadjov2019passivity} and in discrete-time \cite{koshal2016distributed,salehisadaghiani2016distributed,salehisadaghiani2019distributed,tatarenko2019geometric}, where the algorithms were designed based on gradient descent and consensus algorithms. The algorithms in \cite{salehisadaghiani2019distributed} using fixed stepsize schemes may have a faster convergence rate than those in \cite{koshal2016distributed} equipped with vanishing stepsizes. As an NE of a convex game can be equivalently expressed as a zero point of a monotone operator, the authors of \cite{tatarenko2018geometric,pavel2020distributed} proposed distributed algorithms for solving the NE seeking problem by the operator theoretic theory.

All the references mentioned above considered offline games, where both the cost functions and constraints are time-invariant. However, dynamic environments always exist in a multitude of practical applications, such as allocating radio resources and online auction. In this scenario, cost functions and constraints in a game are time-varying and their values and gradients can be accessible only after decisions are made at the current time. These motivate researchers to find distributed online or learning algorithms for seeking NEs or generalized NEs (GNEs). Along this line, the authors of \cite{lu2020online} applied a primal-dual strategy and consensus algorithms to devise a no-regret online algorithm for seeking GNE of a time-varying game, where cost functions are time-varying while nonlinear constraints are invariant.


In this paper, a time-varying noncooperative game is considered, where each player selfishly minimizes its own time-varying cost function subject to time-varying coupled constraints and local feasible set constraints. Individual players only have access to the local cost functions and local constraints, and can communicate with some other players via a connected graph. To solve this problem, we present a distributed online algorithm for seeking GNE based on mirror descent and primal-dual algorithms when gradients of cost functions and nonlinear constraint functions are available, and rigorously prove that Algorithm \ref{alg1} can achieve that dynamic regrets and constraint violation grow sublinearly by appropriately choosing some decreasing stepsizes.
The main contributions of this paper are twofold.
\begin{itemize}
\item[1)] A distributed online algorithm (Algorithm \ref{alg1}) is presented for seeking GNE of a noncooperative game with time-varying cost functions and shared nonlinear constraints, while algorithms in  \cite{facchinei2010generalized}--\!\!\cite{pavel2020distributed} are only applicable to static games and the online algorithm in \cite{lu2020online} is for time-varying games but with time-invariant constraints.
\item[2)] Mirror descent is utilized in designing the algorithm of this paper, which is more applicable and general than projection-based algorithms \cite{lu2020online}, as Bregman divergence is employed.
\end{itemize}

The rest of this paper is structured as follows. In Section \ref{section2}, some preliminaries and the problem formulation, are introduced. Section \ref{section3} presents a distributed online algorithm for seeking GNE.
A numerical example is given to show the effectiveness of the proposed algorithm in Section \ref{section5}. Section \ref{section6} makes a brief conclusion.

{\em Notations.} The symbols $\mathbb{R}$, $\mathbb{R}^m$, and $\mathbb{R}^{n\times m}$ represent the sets of real numbers, $m$-dimensional real column vectors, and $n\times m$ real matrices, respectively. Let $\mathbb{R}^m_+$ be the set of $m$-dimensional nonnegative vectors. $\mathbb{S}^m$ denotes the sphere centered at the origin in $\mathbb{R}^m$, while $\mathbb{B}^m$ is the unit ball centered at the origin in $\mathbb{R}^m$. The symbol ${\bf 1}_m$ (resp. ${\bf 0}_m$) represents an $m$-dimensional vector, whose entries are 1 (resp. 0). For a vector or matrix $A$, the transpose of $A$ is denoted by $A^\top$.  The identity matrix of dimension $m$ is denoted by $I_m$. For an integer $m>0$, let $[m]:=\{1,2,\ldots,m\}$. Let $col(y_1,\ldots,y_m):=(y_1^{\top},\ldots,y_m^{\top})^{\top}$. $P\otimes Q$ denotes the Kronecker product of matrices $P$ and $Q$. For a vector $v\in\mathbb{R}^m$, $[v]_+$ is the projection of $v$ onto $\mathbb{R}^m_+$. $\langle x,y\rangle$ is the standard inner product of $x\in\mathbb{R}^m$ and $y\in\mathbb{R}^m$. For two vectors/matrices $x,y\in\mathbb{R}^m$, the symbol $x\leq y$ means that each entry of $x-y$ is nonpositive, while for two real symmetric matrices $W,P\in\mathbb{R}^{m\times m}$, $W\succeq P$ and $W\succ P$ mean that $W-P$ is positive semi-definite and positive definite, respectively. Given two functions $h_1(\cdot)$ and $h_2(\cdot)$, the notations $h_1=\mathcal{O}(h_2)$ and $h_1={\bf o}(T)$ mean that there exists a positive constant $C>0$ such that $|h_1(x)|\leq Ch_2(x)$ and $\lim\limits_{T\to\infty}\frac{h_1}{T}=0$ for any $x$ in the domain, respectively.


\section{Preliminaries}\label{section2}
\subsection{Graph Theory}\label{subII.A}
 Let ${\mathcal G}=({\mathcal V},{\mathcal E},A)$ be an undirected graph, whose vertex set is ${\mathcal V}=[N]$, the edge set is ${\mathcal E}\subseteq{\mathcal V}\times{\mathcal V}$ and the weighted adjacency matrix is $A=(a_{ij})\in\mathbb{R}^{N\times N}$. For any $i,j\in[N]$ and $i\neq j$, $a_{ij}>0$ if $(j,i)\in{\mathcal E}$ and $a_{ij}=0$ otherwise. In this paper, it is assumed that $a_{ii}>0$ for all $i\in[N]$. $j$ is called a neighbor of $i$ if $(j,i)\in\mathcal{E}$. Denote ${\mathcal N}_i:=\{j:~(j,i)\in{\mathcal E}\}$. A path from node $i_1$ to node $i_l$ is composed of a sequence of edges $(i_h,i_{h+1})$, $h=1,2,\ldots,l-1$. It is said that an undirected graph ${\mathcal G}$ is connected if there exists a path from node $i$ to node $j$ for any vertices $i,j$.

 For communication graph ${\mathcal G}$, the following standard assumptions are imposed in this paper.
 \begin{assumption}\label{assump1}
The undirected graph ${\mathcal G}$ is connected and the adjacency matrix $A$ satisfies that $A^{\top}=A$ and $A{\bf1}_N={\bf1}_N$.
 \end{assumption}

 Let $A_i^-$ be a matrix in $\mathbb{R}^{(N-1)\times(N-1)}$, obtained by deleting the $i$th row and $i$th column of $A$. By Assumption \ref{assump1} and Lemma 3 in \cite{hong2006tracking}, one has that $I_{N-1}-A_i^-\succ0$. Therefore, all the eigenvalues of $A_i^-$ are less than 1. By the Gershgorin circle theorem, it can be easily derived that $\lambda(A_i^-)>-1$. Thus, $-1<\lambda(A_i^-)<1$ for every $i\in[N]$. Denote $\sigma:=\max_{i\in[N]}|\lambda(A_i^-)|\in(0,1)$. Moreover, denote $\sigma_m:=\lambda_{\max}(A-{\bf1}_N{\bf1}_N^{\top}/N),$ where $\lambda_{\max}(A)$ represents the maximum eigenvalue of matrix $A$. Then $0<\sigma_m<1$ \cite{li2020}.
\subsection{Bregman Divergence}
For each player $i\in[N]$, the Bregman divergence $D_{\phi_i}(\xi,\zeta)$ of two points $\xi,\zeta\in\Omega_i\subseteq\mathbb{R}^{n_i}$ is defined as \cite{bregman1967relaxation}
\begin{align}
D_{\phi_i}(\xi,\zeta):=\phi_i(\xi)-\phi_i(\zeta)-\langle\nabla\phi_i(\zeta),\xi-\zeta\rangle,
\end{align} where $\phi_i:\Omega_i\to\mathbb{R}$ is differentiable and $\mu_i$-strongly convex for some constant $\mu_i>0$, i.e., $$\phi_i(\xi)\geq\phi_i(\zeta)+\langle\nabla\phi_i(\zeta),\xi-\zeta\rangle+\frac{\mu_i}{2}\|\xi-\zeta\|^2.$$ Thus, it can be easily derived that $D_{\phi_i}(\cdot,\zeta)$ is $\mu_0$-strongly convex for $\mu_0:=\min\{\mu_1,\ldots,\mu_N\}$, i.e.,
\begin{align}\label{equ2}
D_{\phi_i}(\xi,\zeta)\geq\frac{\mu_0}{2}\|\xi-\zeta\|^2,
\end{align}  and the generalized triangle inequality is satisfied, i.e.,
\begin{align}
&\langle \xi-\zeta, \nabla\phi_i(\zeta)-\nabla\phi_i(\theta)\rangle\nonumber\\
&=D_{\phi_i}(\xi,\theta)-D_{\phi_i}(\xi,\zeta)-D_{\phi_i}(\zeta,\theta).\label{equ3}
\end{align}

Two typical examples of Bregman divergence are the Euclidean distance $D_{\phi_i}(\xi,\zeta)=\|\xi-\zeta\|^2$ generated by $\phi_i(\xi)=\|\xi\|^2$ and the generalized Kullback-Leibler divergence $D_{\phi_i}(\xi,\zeta)=\sum_{j=1}^{n_i}\xi_j\log\frac{\xi_j}{\zeta_j}-\sum_{j=1}^{n_i}\xi_j+\sum_{j=1}^{n_i}\zeta_j$ generated by $\phi_i(\xi)=\sum_{j=1}^{n_i}\xi_j\log \xi_j$. Two mild assumptions on Bregman divergence are given as follows.
\begin{assumption}\label{assump2}
For any $i\in[N]$ and $\xi,\zeta\in\Omega_{i}$, $D_{\phi_i}(\xi,\zeta)$ is Lipschitz with respect to the first variable $\xi\in\Omega_i$, i.e., one can find a positive number $K$ such that for any $\xi_1,\xi_2\in\Omega_{i}$,
\begin{align}
|D_{\phi_i}(\xi_1,\zeta)-D_{\phi_i}(\xi_2,\zeta)|\leq K\|\xi_1-\xi_2\|.
\end{align}
\end{assumption}
\begin{assumption}\label{assump3}
For any $i\in[N]$ and $\xi\in\Omega_{i}$, $D_{\phi_i}(\xi,\cdot):\Omega_{i}\to\mathbb{R}$ is convex, i.e., for any $a\in[0,1]$,
\begin{align}
&D_{\phi_i}(\xi,a\zeta_1+(1-a)\zeta_2)\nonumber\\
&\leq aD_{\phi_i}(\xi,\zeta_1)+(1-a)D_{\phi_i}(\xi,\zeta_2),~\forall \zeta_1,\zeta_2\in\Omega_{i}.
\end{align}
\end{assumption}

Assumption \ref{assump2} is satisfied when $\phi_i(\xi)$ is Lipschitz on $\Omega_{i}$ and also implies that for any $\xi,\zeta\in\Omega_i$, $D_{\phi_i}(\xi,\zeta)=|D_{\phi_i}(\xi,\zeta)-D_{\phi_i}(\zeta,\zeta)|\leq K\|\xi-\zeta\|$. Assumption \ref{assump3} is crucial to derive the main results in this paper and a sufficient condition given in \cite{bauschke2001joint} for guaranteeing Assumption \ref{assump3} is that $\phi_i(\xi)$ is thrice continuously differentiable and satisfies
$H_{\phi_i}(\xi)\succeq 0$, $H_{\phi_{i}}(\xi)+\nabla H_{\phi_i}(\xi)(\xi-\zeta)\succeq0$, where $\xi,\zeta\in\Omega_{i}$ and $H_{\phi_i}$ represents the Hessian matrix of $\phi_i$.
\subsection{Problem Formulation}
Denote by $\Gamma(\mathcal{V},\Omega,J)$ a noncooperative game with $N$ players, where $\mathcal{V}:=[N]$ is the set of players, $\Omega:=\Omega_1\times\cdots\times\Omega_N$ represents the strategy set of players with $\Omega_i\subseteq\mathbb{R}^{n_i}$ being the private action set of player $i$, and $J=(J_1,\ldots,J_N)$ is the cost function with $J_i$ being the cost function of player $i$. Denote by $x=col(x_1,\ldots,x_N)$ the joint action, where $x_i$ is the action of player $i$, $i\in[N]$. Denote by $x_{-i}:=col(x_1,\ldots,x_{i-1},x_{i+1},\ldots,x_N)$ the joint action of all the players except $i$. For a game $\Gamma(\mathcal{V},\Omega,J)$, a strategy profile $x^*=(x_1^*,\ldots,x_N^*)\in\Omega$ is called an NE if for any $i\in[N]$, there holds that
\begin{align}
J_i(x_i^*,x_{-i}^*)\leq J_i(x_i,x_{-i}^*), ~\forall x_i\in\Omega_i.
\end{align}
Moreover, if $\Omega_i$ depends on other players' actions, then the NE $x^*$ is termed a GNE.

In this paper, a time-varying game $\Gamma(\mathcal{V},\Omega_t,J_t)$ under dynamic environments is studied and the players can receive their neighbors' information following a fixed graph ${\cal G}=(\mathcal{V},{\mathcal E},A)$, satisfying Assumption \ref{assump1}. The cost function $J_{t}=(J_{1,t},\ldots,J_{N,t})$ and the action set $\Omega_t=\Omega_{0,t}\bigcap(\Omega_{1}\times\cdots\times\Omega_{N})$ are time-varying, where $\Omega_{0,t}$ is the shared convex constraint $\Omega_{0,t}:=\{x\in\mathbb{R}^n\mid g_t(x_t):=\sum_{i=1}^Ng_{i,t}(x_i)\leq{\bf 0}_m\}$ and $\Omega_{i}$ is the private action set constraint of player $i\in[N]$. Here, $n:=\sum_{i=1}^Nn_i$ and $g_{i,t}:\mathbb{R}^{n_i}\to\mathbb{R}^m$.
At each time $t$, each player $i$ has the purpose to solve the following optimization problem:
\begin{align}
&\min\limits_{x_{i,t}\in\Omega_{i}}~~~~ J_{i,t}(x_{i,t},x_{-i,t})\nonumber\\
&\text{subject to} ~ ~~g_t(x_t)\leq {\bf 0}_m.
\end{align}

For each player $i\in[N]$, denote $g_{i,t}:=col(g_{i1,t},\ldots,g_{im,t})$, where $g_{ij,t}:\mathbb{R}^{n_i}\to\mathbb{R}$ for $j\in[m]$. Some standard assumptions are needed, which are also made in \cite{salehisadaghiani2016distributed,lu2020online}.
\begin{assumption}\label{assump4}
For each $i\in[N]$, the non-empty set $\Omega_i$ is compact and convex. The differentiable function $J_{i,t}(\cdot,x_{-i}):\Omega_i\to\mathbb{R}$ is convex for any $x_{-i}\in\mathbb{R}^{n-n_i}$. $g_{ij,t}(x_i)$ is differentiable and convex for any $x_i\in\mathbb{R}^{n_i}$. Moreover, the constraint set $\Omega_t$ is assumed to be non-empty and Slater's constraint qualification is satisfied.
\end{assumption}

Under Assumption \ref{assump4}, $\|x_i\|$, $\|J_{i,t}(x)\|$ and $\|g_{i,t}(x_i)\|$ are bounded for any $i\in[N]$, $x_i\in\Omega_{i}$ and $x_{-i}\in\mathbb{R}^{n-n_i}$. Thus, it can be assumed that $L>0,M>0$ can be found such that
\begin{align}
\|x_i\|\leq L, \|J_{i,t}(x)\|\leq L, \|g_{i,t}(x_i)\|\leq L,\label{equ11}\\
\|\nabla_iJ_{i,t}(x)\|\leq M,~\|\nabla g_{i,t}(x_i)\|\leq M,\label{equ12}
\end{align} where $\nabla g_{i,t}(x_i):=col((\nabla g_{i1}(x_i))^{\top},\ldots,(\nabla g_{im}(x_i))^{\top})$ and $\nabla_iJ_{i,t}(x_i,x_{-i}):=\frac{\partial J_{i,t}(x_i,x_{-i})}{\partial x_i}$.

Define
$
F_t(x):=col(\nabla_1J_{1,t}(x),\ldots,\nabla_NJ_{N,t}(x)),
$ which is called the {\em pseudo-gradient mapping} of the game $\Gamma(\mathcal{V},\Omega_t,J_t)$.
\begin{assumption}\label{assump5}
For any $i\in[N]$, $x=(x_i,x_{-i})$ and $y=(y_i,y_{-i})$, where $x_i,y_i\in\Omega_{i}$ and $x_{-i},y_{-i}\in\mathbb{R}^{n-n_i}$, the mapping $\nabla_iJ_{i,t}(x_{i},x_{-i})$ is $H$-Lipschitz continuous, i.e.,
\begin{align}
\|\nabla_iJ_{i,t}(x_{i},x_{-i})-\nabla_iJ_{i,t}(y_{i},y_{-i})\|\leq H\|x-y\|.
\end{align}
\end{assumption}
\begin{assumption}\label{assump6}
The mapping $F_t(x)$ is $\mu$-strongly monotone on the set $\Omega_{1}\times\cdots\times\Omega_{N}$ for a constant $\mu>0$, i.e., for any $x,y\in\Omega_{1}\times\cdots\times\Omega_{N}$,
\begin{align}
(F_t(x)-F_t(y))^{\top}(x-y)\geq\mu\|x-y\|^2.
\end{align}
\end{assumption}

Note that the cost functions are convex and differentiable. It can be obtained from Theorem 3.9 in \cite{facchinei2009nash} that at any time $t$, a solution to the following variational inequality:
\begin{align}\label{equ8}
(F_t(x^*_t))^{\top}(x-x^*_t)\geq 0, ~\text{for all}~x\in\Omega_t,
\end{align}
is a GNE of $\Gamma(\mathcal{V},\Omega_t,J_t)$, and this GNE $x^*_t$ is also called a variational GNE. In addition, we have that the inequality (\ref{equ8}) has a unique solution by Assumption \ref{assump6}. Therefore, the existence and uniqueness of the variational GNE can be guaranteed by Assumptions \ref{assump4} and \ref{assump6}. It is noted that finding all GNEs is very difficult even if the game is offline. Accordingly, we will discuss on tracking the unique variational GNE as done in \cite{liang2017distributed,pavel2020distributed} since the unique variational GNE enjoys good stability and has no price discrimination from the perspective of economics.

The objective of this paper is to devise distributed online algorithms to mimic the performance of its offline counterpart. By the definition of GNEs, $x_t^*=(x_{i,t}^*,x_{-i,t}^*)$ is the GNE of the time-varying game $\Gamma(\mathcal{V},\Omega_t,J_t)$ if and only if $x_{i,t}^*$ is a solution to the following optimization problem
\begin{align}
&\min\limits_{x_{i,t}\in\Omega_{i}}~~~~ J_{i,t}(x_{i,t},x^*_{-i,t})\nonumber\\
&\text{subject to} ~~~g_{i,t}(x_{i,t})+\sum_{j=1,j\neq i}^Ng_{j,t}(x^*_{j,t})\leq {\bf 0}_m.\label{equ14}
\end{align} The {\em dynamic regret} of player $i\in[N]$ is defined as
 \begin{align}\label{equ10}
 Reg_i(T):=\sum_{t=1}^T(J_{i,t}(x_{i,t},x^*_{-i,t})-J_{i,t}(x^*_{i,t},x^*_{-i,t})),
 \end{align} where $x^*_t=(x_{i,t}^*,x_{-i,t}^*)$ is the variational GNE of $\Gamma(\mathcal{V},\Omega_t,J_t)$ at time $t$ and $T$ is the learning time. For a decision sequence $\{x_{1},\ldots,x_{T}\}$, the constraint violation measure is given as
 \begin{align}\label{e11}
 R_g(T):=\left\|\left[\sum\limits_{t=1}^Tg_t(x_t)\right]_+\right\|.
 \end{align}

It is said that an online algorithm is announced ``good'' or is a no-regret algorithm if all the regrets of players in (\ref{equ10}) and the accumulation of constraint violations in (\ref{e11}) increase sublinearly, i.e.,
$
Reg_i(T)={\bf o}(T),~R_g(T)={\bf o}(T).
$

However, it could be impossible to track the GNE if the variational GNE sequence $\{x^*_1,\ldots,x^*_T\}$ of the considered time-varying game fluctuates drastically. For convenience, motivated by \cite{hall2015online,lu2020online}, the following accumulation is adopted:
\begin{align}
\Phi_T^*:=\sum_{t=1}^T\|x^*_{t+1}-x^*_t\|.
\end{align}

\section{Main Result}\label{section3}
For the time-varying game $\Gamma(\mathcal{V},\Omega_t,J_t)$, each player $i$ only has the information from its neighbors and the information associated with $J_{i,t}$, $\Omega_i$ and $g_{i,t}$. In this setting, a distributed online algorithm based on mirror descent will be presented.

For each $i\in[N]$, define an augmented Lagrangian function at time $t$ as
$
{\cal L}_{i,t}(x_{i,t},\lambda_t;x_{-i,t}):=J_{i,t}(x_{i,t},x_{-i,t})+\gamma_t\lambda^{\top}_tg_t(x_t)-\frac{\beta_{t}}{2}\|\lambda_t\|^2,
$ where $\lambda_t\in\mathbb{R}^m_+$ is Lagrange multiplier or the dual variable, $\gamma_t>0$ is a stepsize, and $\beta_t>0$ is the regularization parameter. Invoking Lemma 1 in \cite{nedic2009approximate}, the optimal dual variable $\lambda_t^*$ is bounded, that is, there exists a positive constant $\Lambda>0$ such that
\begin{align}
\|\lambda_t^*\|\leq\Lambda.\label{e20}
\end{align} Inspired by the dynamic mirror descent for online optimization \cite{hall2015online}, to track the GNE online, a primal-dual dynamic mirror descent under full-decision information, i.e., each player has access to the others' decisions $x_{-i}$, can be designed as
\begin{align}
x_{i,t+1}&=\arg\min\limits_{x\in\Omega_i}\{\alpha_{t}\langle x,\nabla_iJ_{i,t}(x_{i,t},x_{-i,t})\rangle\nonumber\\
&~~~~~~~+\alpha_{t}\langle x,\gamma_t(\nabla g_{i,t}(x_{i,t}))^{\top}\lambda_t\rangle+D_{\phi_i}(x,x_{i,t})\},\label{equ20}\\
\lambda_{t+1}&=\left[\lambda_t+\gamma_t(\gamma_tg_t(x_t)-\beta_t\lambda_t)\right]_+,\label{equ21}
\end{align} where $\alpha_t>0$ and $0<\gamma_t<1$ ($\alpha_0=\gamma_0=1$) are time-varying stepsizes utilized in the primal and dual iterations. The main drawbacks of this algorithm are that each player should know all the others' decisions, the common Lagrange multiplier $\lambda_t$ and the nonlinear constraint function $g_t$. Moreover, the stepsize sequences are designed based on the upper bounds of the cost functions and constraint functions, along with their subgradients. In order to avoid such disadvantages, at time slot $t$, let $x_{ij,t}$ and $\lambda_{i,t}$ be estimates of the strategy of player $j$ and the global Lagrange multiplier by player $i$, respectively. $x_{ii,t}:=x_{i,t}$. Motivated by the algorithm proposed in \cite{neely2017online}, by modifying (\ref{equ20}) and (\ref{equ21}), a distributed online primal-dual dynamic mirror descent algorithm as in Algorithm \ref{alg1} is designed to learn the variational GNE of the time-varying game $\Gamma(\mathcal{V},\Omega_t,J_t)$ under partial-decision information.
In order to execute Algorithm \ref{alg1}, at each time slot $t$, every player $i$ needs to know $\nabla_iJ_{i,t}({\bf x}_{i,t})$, $g_{i,t}(x_{i,t})$ and $\nabla g_{i,t}(x_{i,t})$ rather than the full information of $J_{i,t}$ and $g_{i,t}$, which is similar to most online algorithms for optimization and games \cite{mahdavi2012trading,jenatton2016adaptive,yuan2018online,lu2020online,li2020distributed,li2020distributed_online}.

\begin{algorithm}[!htbp]\caption{Distributed Online Primal-Dual Dynamic Mirror Descent}\label{alg1}
Each player $i$ maintains vector variables $x_{ih,t}\in\mathbb{R}^{n_h}$ and $\lambda_{i,t}\in\mathbb{R}^{m}$ at iteration $t\in[T]$.

 {\bf Initialization:} For any $i\in[N]$, initialize $x_{i,1}\in\Omega_i$ arbitrarily, $x_{ih,1}={\bf0}_{n_h}(h\neq i)$, and $\lambda_{i,1}={\bf0}_m$.

{\bf Iteration:} For $t\geq 1$, every player $i$ processes the following update:
\begin{align}
x_{ih,t+1}&=\sum\limits_{k=1}^Na_{ik}x_{kh,t},~~h\neq i,\\
\tilde{x}_{i,t+1}&=\arg\min\limits_{x\in\Omega_i}\{\alpha_t\langle x,\nabla_i J_{i,t}({\bf x}_{i,t})+\gamma_t(\nabla g_{i,t}(x_{i,t}))^{\top}\tilde{\lambda}_{i,t}\rangle\nonumber\\
&~~~~~~~~~~~~~~+D_{\phi_i}(x,x_{i,t})\},\label{equ23}\\
x_{i,t+1}&=(1-\alpha_t)x_{i,t}+\alpha_t\tilde{x}_{i,t+1},\label{equ24}\\
\lambda_{i,t+1}&=\left[\tilde{\lambda}_{i,t}+\gamma_t(\gamma_tb_{i,t+1}-\beta_t\tilde{\lambda}_{i,t})\right]_+,\label{equ25}
\end{align} where ${\bf x}_{i,t}:=col(x_{i1,t},\ldots,x_{iN,t})$, $\tilde{\lambda}_{i,t}:=\sum_{j=1}^Na_{ij}\lambda_{j,t}$, $b_{i,t+1}:=\nabla g_{i,t}(x_{i,t})(\tilde{x}_{i,t+1}-x_{i,t})+g_{i,t}(x_{i,t})$, $a_{ij}$ is the $(i,j)$th element of $A$, and $\alpha_t,\beta_t,\gamma_t$, satisfying  $\alpha_0=\beta_0=\gamma_0=1$ and $\frac{\gamma_t}{\beta_t}\leq\frac{\gamma_{t+1}}{\beta_{t+1}}$, are the stepsizes to be determined.
\end{algorithm}

It is noted that Algorithm \ref{alg1} subsumes a few interesting frameworks. For example, let $\phi_i(\xi)=\frac{1}{2}\|\xi\|^2$, then the associated Bregman divergence is $D_{\phi_i}(\xi,\zeta)=\frac{1}{2}\|\xi-\zeta\|^2$. In this setting, Algorithm \ref{alg1} reduces to a variant of the distributed online projection-based algorithm in \cite{lu2020online}. Assume that the constraint set is $\Omega_i=\{\xi\in\mathbb{R}^{n_i}: \sum_{i=1}^{n_i}[\xi]_i=1, \xi\geq{\bf 0}_{n_i}\}$, where $[\xi]_i$ represents the $i$th element of $\xi$. Let $\phi_i(\xi)=\sum_{i=1}^{n_i}[\xi]_i\ln[\xi]_i$, then the associated Bregman divergence is the Kullback-Leibler divergence $D_{\phi_i}(\xi,\zeta)=\sum_{i=1}^{n_i}[\xi]_i\ln\frac{[\xi]_i}{[\zeta]_i}$. In this case, Algorithm \ref{alg1} reduces to a revised version of the distributed entropic descent algorithm in \cite{beck2003mirror}. However, the projection step in the projection-based algorithm is not easy to be solved explicitly.

In what follows, some necessary lemmas are presented. First, a result on the bounds of the estimate errors of players' strategies is given as follows.
\begin{lemma}\label{lemma1}
If Assumptions \ref{assump1} and \ref{assump4} are satisfied, then, for any $i\in[N]$ and $t\in[T]$,
\begin{align}
\|e_{i,t}\|\leq2\sqrt{N-1}L\sum\limits_{s=0}^{t-1}\sigma^s\alpha_{t-s-1},\label{e27}
\end{align} where $e_{i,t}:=col(e_{1i,t},\ldots,e_{(i-1)i,t},e_{(i+1)i,t},\ldots,e_{Ni,t})$, $e_{ih,t}:=x_{ih,t}-x_{h,t}$ and $0<\sigma<1$ is defined in Subsection \ref{subII.A}.
\end{lemma}

{\em Proof:}
See Appendix \ref{A0}. \hfill$\blacksquare$
\begin{lemma}\label{lemma2}
If Assumptions \ref{assump1} and \ref{assump4} hold, then for any $i\in[N]$ and $t\in[T]$, $\lambda_{i,t}$ and $\tilde{\lambda}_{i,t}$ generated by Algorithm \ref{alg1} satisfy
\begin{align}
\|\lambda_{i,t}\|&\leq\frac{L\gamma_t}{\beta_t},\label{equ27}\\
\|\tilde{\lambda}_{i,t}\|&\leq\frac{L\gamma_t}{\beta_t},\label{equ28}\\
\|\tilde{\lambda}_{i,t}-\overline{\lambda}_t\|&\leq\sqrt{N}F\sum\limits_{s=0}^{t-1}\sigma_m^{s}\gamma_{t-1-s}^2,\label{equ29}\\
\frac{\Xi_{t+1}}{2\gamma_t}&\leq \frac{NF^2}{2}\gamma_t^3+\gamma_t(\overline{\lambda}_t-\lambda)^{\top}g_t(x_t)+c_1(t)\nonumber\\
&~~~+c_2(t)+N\left(\frac{M^2\alpha_t\gamma_t^2}{\mu_0}+\frac{\beta_t}{2}\right)\|\lambda\|^2,\label{equ30}
\end{align} where $\lambda\in\mathbb{R}^m_+$, $\overline{\lambda}_t:=\frac{1}{N}\sum_{i=1}^N\lambda_{i,t}$, $F:=2LM+2L$, $c_1(t):=N\sqrt{N}FL\gamma_t\sum\limits_{s=0}^{t-1}\sigma_m^{s}\gamma_{t-1-s}^2$, and
\begin{align}
\Xi_{t+1}&:=\sum_{i=1}^N\left[\|\lambda_{i,t+1}-\lambda\|^2-(1-\beta_t\gamma_t)\|\lambda_{i,t}-\lambda\|^2\right],\nonumber\\
c_2(t)&:=\gamma_t\sum_{i=1}^N(\tilde{\lambda}_{i,t})^{\top}\nabla g_{i,t}(x_{i,t})(\tilde{x}_{i,t+1}-x_{i,t})\nonumber\\
&~~~~+\frac{\mu_0}{4\alpha_t}\sum_{i=1}^N\|\tilde{x}_{i,t+1}-x_{i,t}\|^2.
\end{align}
\end{lemma}

{\em Proof:}
See Appendix \ref{A}. \hfill$\blacksquare$
\begin{lemma}\label{lemma3}
Under Assumptions \ref{assump1}--\ref{assump6}, for any $i\in[N]$ and $t\in[T]$, $x_{i,t}$ generated by Algorithm \ref{alg1} satisfies
\begin{align}
&\mu\sum_{t=1}^T\|x_{t}-x_{t}^*\|^2\nonumber\\
&\leq\sum_{t=1}^T\frac{1}{\alpha_t^2}\sum_{i=1}^N\left[D_{\phi_i}(x_{i,t}^*,x_{i,t})-D_{\phi_i}(x_{i,t+1}^*,x_{i,t+1})\right]\nonumber\\
&~~~+\frac{\sqrt{N}K}{\alpha_T^2}\Phi_T^*+C_1\sum_{t=1}^T\gamma_{t-1}+C_2\sum_{t=1}^T\alpha_t+\frac{N}{2}\|\lambda\|^2\nonumber\\
&~~~-\lambda^{\top}\sum_{t=1}^T\gamma_tg_t(x_t)+N\sum_{t=1}^T\left(\frac{M^2}{\mu_0}\alpha_t\gamma_t^2+\frac{1}{2}\beta_t\right)\|\lambda\|^2\nonumber\\
&~~~+\frac{1}{2}\sum\limits_{t=1}^T\left(\frac{1}{\gamma_t}-\frac{1}{\gamma_{t-1}}-\beta_t\right)\sum_{i=1}^N\|\lambda_{i,t}-\lambda\|^2,\label{equ34}
\end{align} where $C_1:=\frac{3N\sqrt{N}FL}{1-\sigma_m}+2NLM\Lambda+\frac{NF^2}{2}$ and $C_2:=\frac{NM^2}{\mu_0}+\frac{4N^2L^2H}{1-\sigma}$.
\end{lemma}

{\em Proof:}
See Appendix \ref{B}. \hfill$\blacksquare$
\begin{lemma}\label{lemma4}
Under Assumptions \ref{assump1}--\ref{assump6}, for any $i\in[N]$, the dynamic regret (\ref{equ10}) generated by Algorithm \ref{alg1} is bounded by
\begin{align}
Reg_i(T)&\leq M\sqrt{TB_1(T)+TB_2(T)+\frac{\sqrt{N}K}{\mu\alpha_T^2}T\Phi_T^*}\label{e35}
\end{align} and the accumulated constraint violation satisfies
\begin{align}
&\left\|\left[\sum\limits_{t=1}^T\gamma_tg_t(x_t)\right]_+\right\|^2\nonumber\\
&\leq4\mu NL^2TB_3(T)+\mu B_1(T)B_3(T)+\frac{\sqrt{N}K}{\alpha_T^2}B_3(T)\Phi_T^*\nonumber\\
&~~~+\frac{B_3(T)}{2}\sum\limits_{t=1}^T\left(\frac{1}{\gamma_t}-\frac{1}{\gamma_{t-1}}-\beta_t\right)\sum_{i=1}^N\|\lambda_{i,t}-\lambda_c\|^2, \label{e36}
\end{align} where $\lambda_c:=\frac{2\left[\sum_{t=1}^T\gamma_tg_t(x_t)\right]_+}{B_3(T)}$ and
\begin{align*}
B_1(T)&:=\frac{2NLK}{\mu\alpha_{T+1}^2}+\frac{C_1}{\mu}\sum_{t=1}^T\gamma_{t-1}+\frac{C_2}{\mu}\sum_{t=1}^T\alpha_t,\nonumber\\
B_2(T)&:=\frac{1}{2\mu}\sum_{t=1}^T\sum\limits_{i=1}^N\left(\frac{1}{\gamma_t}-\frac{1}{\gamma_{t-1}}-\beta_t\right)\|\lambda_{i,t}\|^2\nonumber\\
B_3(T)&:=4N\left(\frac{1}{2}+\sum_{t=1}^T\left(\frac{M^2}{\mu_0}\alpha_t\gamma_t^2+\frac{1}{2}\beta_t\right)\right).
\end{align*}
\end{lemma}

{\em Proof:}
See Appendix \ref{C}. \hfill$\blacksquare$

Then it is ready to give the main result on the bounds of the regrets and constraint violation for Algorithm \ref{alg1}.
\begin{theorem}\label{thm1}
If Assumptions \ref{assump1}--\ref{assump6} hold, then for each $i\in[N]$ and the sequence $\{x_{i,1},\ldots,x_{i,T}\}$ generated by Algorithm \ref{alg1} with
$$\alpha_t=\frac{1}{t^{a_1}},\beta_t=\frac{1}{(T+1)^{a_2}},\gamma_t=\frac{1}{(T+1)^{1-a_2}},$$ where $0<a_1<\frac{1}{2}$ and $\frac{2}{3}<a_2<1$, there hold
\begin{align}
Reg_{i}(T)&={\mathcal O}\left(T^{\max\{\frac{1}{2}+a_1,\frac{1}{2}+\frac{a_2}{2},1-\frac{a_1}{2}\}}\right)\nonumber\\
&~~~+\mathcal{O}\left(T^{\frac{1}{2}+a_1}\sqrt{\Phi_T^*}\right),\label{equ35}\\
R_g(T)&= \mathcal{O}\left(T^{\max\{1-\frac{a_1}{2},2-\frac{3}{2}a_2\}}\right)\nonumber\\
&~~~+\mathcal{O}\left(T^{\max\{\frac{1}{2}+\frac{a_1}{2},\frac{3}{2}+a_1-\frac{3a_2}{2}\}}\sqrt{\Phi_T^*}\right).\label{equ36}
\end{align}
\end{theorem}

{\em Proof:}
For the selected parameters $\beta_t$ and $\gamma_t$, it is easy to verify that $\frac{1}{\gamma_t}-\frac{1}{\gamma_{t-1}}-\beta_t\leq0$ and $\frac{\gamma_t}{\beta_t}\leq\frac{\gamma_{t+1}}{\beta_{t+1}}$.

In addition, for any constant $0<a\neq1$ and positive integer $T$, it holds that
\begin{align}\label{e90}
\sum\limits_{t=1}^T\frac{1}{t^a}\leq1+\int_{1}^T\frac{1}{t^{a}}dt\leq1+\frac{T^{1-a}-1}{1-a}\leq\frac{T^{1-a}}{1-a}.
\end{align} Then,
\begin{align}
\sum_{t=1}^T\alpha_t&\leq\frac{T^{1-a_1}}{1-a_1}.\label{equ58}
\end{align} Thus, one can derive that
\begin{align}
B_1(T)&=\mathcal{O}(T^{2a_1}+T^{a_2}+T^{1-a_1}),\\
B_3(T)&=\mathcal{O}(T^{2a_2-a_1-1}+T^{1-a_2}).
\end{align}

Then one can obtain from Lemma \ref{lemma4} that
\begin{align*}
Reg_i(T)&=\mathcal{O}\left(\sqrt{TB_1(T)+T^{1+2a_1}\Phi_T^*}\right)\\
&= \mathcal{O}\left(\sqrt{T^{\max\{1+2a_1,1+a_2,2-a_1\}}+T^{1+2a_1}\Phi_T^*}\right),
\end{align*} and
\begin{align*}
&\left\|\left[\sum\limits_{t=1}^T\frac{1}{T^{1-a_2}}g_t(x_t)\right]_+\right\|^2\\
&=\mathcal{O}\left(TB_3(T)+B_1(T)B_3(T)+T^{2a_1}B_3(T)\Phi_T^*\right)\\
&=\mathcal{O}\left(T^{\max\{2a_2-a_1,2-a_2\}}+T^{\max\{2a_2+a_1-1,1+2a_1-a_2\}}\Phi_T^*\right)
\end{align*} which implies that (\ref{equ36}) holds.
Then the result is proved. \hfill$\blacksquare$

\begin{corollary}
 Under Assumptions \ref{assump1}--\ref{assump6}, Algorithm \ref{alg1} achieves sublinearly bounded regrets and constraint violation if $\sqrt{\Phi_T^*}$ grows sublinearly with a known order in $[0,1)$.
\end{corollary}

{\em Proof:}
Assume that there exists a known constant $\tau\in[0,1)$ such that $\Phi_T^*=\mathcal{O}({T}^{2\tau})$, then setting $a_1\in(0,\frac{1}{2}-\tau)$ and $a_2\in(\frac{2}{3},1)$ in Theorem \ref{thm1} yields that
$
Reg_i(T)={\bf o}(T),~R_g(T)={\bf o}(T).
$  \hfill$\blacksquare$
\begin{corollary}
 Under Assumptions \ref{assump1}--\ref{assump6}, Algorithm \ref{alg1} achieves optimal bounds for the dynamic regrets and constraint violation as
\begin{align*}
Reg_i(T)&=\mathcal{O}(T^{\frac{7}{8}}+T^{\frac{5}{6}}\sqrt{\Phi_T^*}),\\
R_g(T)&=\mathcal{O}(T^{\frac{7}{8}}+T^{\frac{17}{24}}\sqrt{\Phi_T^*}).
\end{align*}
\end{corollary}

{\em Proof:}
As $\frac{1}{2}+a_1$ and $1-\frac{a_1}{2}$ grow when increasing $a_1$ and decreasing $a_1$, respectively, it is reasonable to set $\frac{1}{2}+a_1=1-\frac{a_1}{2}$, i.e., $a_1=\frac{1}{3}$. Then $\max\{\frac{1}{2}+a_1,1-\frac{a_1}{2}\}$ has the minimal value  $\frac{5}{6}$. It can be easily seen that $\frac{1}{2}+\frac{a_2}{2}>\frac{5}{6}$ due to $\frac{2}{3}<a_2<1$. As for the three terms $\frac{1}{2}+\frac{a_2}{2}$, $1-\frac{a_1}{2}$ and $2-\frac{3}{2}a_2$, by letting $\frac{1}{2}+\frac{a_2}{2}=2-\frac{3}{2}a_2$, i.e., $a_2=\frac{3}{4}$, the minimal value of $\max\{\frac{1}{2}+\frac{a_2}{2}, 1-\frac{a_1}{2}, 2-\frac{3}{2}a_2\}$ is $\frac{7}{8}$. The result is thus proved by Theorem \ref{thm1}. \hfill$\blacksquare$
\begin{remark}
A distributed online algorithm was presented to seek the variational GNE with time-invariant nonlinear constraints in \cite{lu2020online}, while it is not applicable to the case of time-varying constraints. The presented algorithm in this paper can deal with the time-varying constraints based on mirror descent, however leads to the difficulty in obtained the bounds of the accumulated errors between $x_t$ and $x_t^*$, and then results in the complexity in the derivation.

\end{remark}
\begin{remark}
In Algorithm \ref{alg1}, the time-varying parameters $\beta_t$ and $\gamma_t$ are designed based on the learning time $T$, which requires that all the players must know $T$ in advance. However, when $T$ is unknown, the classical ``doubling trick'' technique can be leveraged to achieve the similar result, as done in \cite{yu2020low}.
\end{remark}

\section{A Numerical Example}\label{section5}
\begin{figure}[!ht]
 \centering
  \includegraphics[width=2.5in]{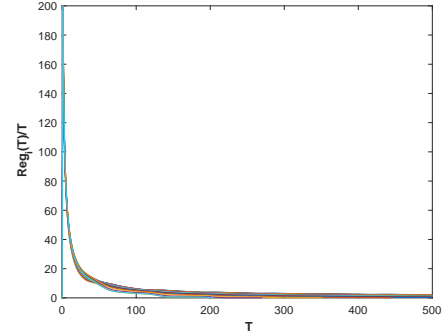}
  \caption{The trajectories of the average regrets $Reg_i(T)/T$, $i\in[N]$ by Algorithm \ref{alg1}.}\label{fig1}
\end{figure}

\begin{figure}[!ht]
 \centering
  \includegraphics[width=2.5in]{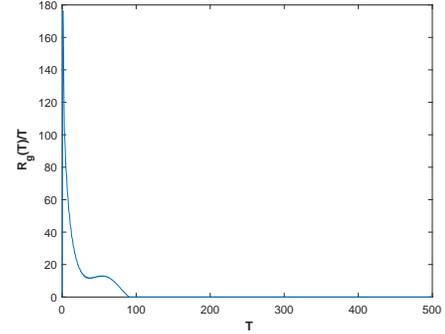}
  \caption{The trajectories of the average violation $R_g(T)/T$, $i\in[N]$ by Algorithm \ref{alg1}.}\label{fig2}
\end{figure}

In this section, a time-varying Nash-Cournot game $\Gamma(\mathcal{V},\Omega_t,J_t)$ with production constraints and market capacity constraints is used to illustrate the feasibility of the obtained algorithm. Similar to \cite{lu2020online}, we consider a Nash-Cournot game, in which there are $N=20$ firms communicating with each other via a connected graph $\mathcal{G}$. Denote by $x_{i,t}\in\mathbb{R}$ the quality produced by firm $i$ at time $t$. In view of some uncertain and changeable factors such as marginal costs and demand for orders, the demand cost and the production cost may be time-varying. Assume that the production cost and the demand price of firm $i$ are $p_{i,t}(x_{i,t})=x_{i,t}(\sin(t/12)+1)$ and $d_{i,t}(x_{t})=22+i/9-0.5i\sin(t/12)-\sum_{j=1}^Nx_{j,t}$, respectively. Then, the overall cost function of firm $i$ is $J_{i,t}(x_{i,t},x_{-i,t})=p_{i,t}(x_{i,t})-x_{i,t}d_{i,t}(x_{t})$ for $i\in[N]$ and $t\in[T]$. In addtion, the production quality constraint of firm $i$ is $x_{i,t}\in\Omega_i:=[0,30]$, while the market capacity constraint is the coupled inequality constraint $\sum_{i=1}^Nx_{i,t}\leq\sum_{i=1}^Nb_{i,t}$, where $b_{i,t}=10+\sin(t/12)$ is the local bound available to firm $i$. In the offline and centralized setting, the GNE can be calculated as $x_{i,t}^*=P_{\Omega_i}(\xi_{i,t})$, where $\xi_{i,t}:=\frac{1}{9}(i-10)+\frac{1}{2}(10-i)\sin\frac{t}{12}$. In the online and distributed setting, by Algorithm \ref{alg1}, set initial states $x_{i,0}\in\Omega_{i}$ randomly, $x_{ih,1}={\bf0}_{n_h}(h\neq i)$, and $\lambda_{i,1}={\bf0}_m$. Choose $a_1=0.2$ and $a_2=0.8$. $Reg_i(T)/T$, $i\in[N]$, and $R_g(T)/T$ are shown in Fig. \ref{fig1} and Fig. \ref{fig2}, respectively.


From these figures, one can see that the average regret $Reg_i(T)/T$ and the average violation $R_g(T)/T$ decay to zero as iteration goes on. That is, the regrets $Reg_i(T)$, $i\in[N]$, and the violation $R_g(T)$ increase sublinearly, which are consistent with Theorem \ref{thm1}.

\section{Conclusion}\label{section6}
In this paper, distributed GNE seeking for noncooperative games with time-varying cost functions and time-varying general convex constraints was investigated. A novel distributed online algorithm was devised based on mirror descent and a primal-dual strategy. It was rigorously proved that the presented algorithm could achieve sublinearly bounded regrets and constraint violation by appropriately choosing decreasing stepsizes. Future research of interest is to develop new online algorithms to improve the bounds of dynamic regrets and constraint violation.

\begin{appendix}
\subsection{Proof of Lemma \ref{lemma1}}\label{A0}
By Lemma 2 in \cite{lu2020online}, one has that
\begin{align*}
\|e_{i,t}\|\leq\sigma^{t-1}\|e_{i,1}\|+2\sqrt{N-1}L\sum\limits_{s=0}^{t-2}\sigma^s\alpha_{t-s-1}.
\end{align*} Since $e_{hi,1}=x_{hi,1}-x_{i,1}=-x_{i,1}$, the norm $\|e_{i,1}\|$ satisfies
\begin{align*}
\|e_{i,1}\|=\sqrt{N-1}\|x_{i,1}\|\leq\sqrt{N-1}L.
\end{align*} Then Lemma \ref{lemma1} can be proved. \hfill$\blacksquare$
\subsection{Proof of Lemma \ref{lemma2}}\label{A}
It can be easily proved that (\ref{equ27}), (\ref{equ28}) and (\ref{equ30}) hold following the proof of Lemma 2 in \cite{yi2020distributed}. Therefore, it suffices to prove (\ref{equ29}). From the i
teration (\ref{equ25}), one has
\begin{align}
\lambda_{i,t+1}=\sum_{j=1}^Na_{ij}\lambda_{j,t}+\epsilon_{i,t},\end{align}
where $\epsilon_{i,t}:=\left[\tilde{\lambda}_{i,t}+\gamma_t(\gamma_tb_{i,t+1}-\beta_t\tilde{\lambda}_{i,t})\right]_+-\tilde{\lambda}_{i,t}$, $i\in[N]$. Denote $\lambda_t:=col(\lambda_{1,t},\ldots,\lambda_{N,t})$ and $\epsilon_t=col(\epsilon_{1,t},\ldots,\epsilon_{N,t})$. Then
\begin{align}
\lambda_{t+1}-{\bf1}_N\otimes\overline{\lambda}_{t+1}&=((A-\frac{1}{N}{\bf1}_N{\bf1}_N^{\top})\otimes I_m)(\lambda_{t}-{\bf1}_N\otimes\overline{\lambda}_t)\nonumber\\
&~~~+((I-\frac{1}{N}{\bf1}_N{\bf1}_N^{\top})\otimes I_m)\epsilon_t. \label{equ40}
\end{align} Under Assumption \ref{assump1}, one has $0<\sigma_m<1$, where $\sigma_m=\lambda_{\max}(A-\frac{1}{N}{\bf1}_N{\bf1}_N^{\top})$. Consequently, taking norm on both sides of (\ref{equ40}) yields
\begin{align}
\|\lambda_{t+1}-{\bf1}_N\otimes\overline{\lambda}_{t+1}\|
&\leq\sigma_m\|\lambda_{t}-{\bf1}_N\otimes\overline{\lambda}_t\|+\|\epsilon_t\|.
\end{align}  Note that
\begin{align}
\|\epsilon_{i,t}\|&= \left\|\left[\tilde{\lambda}_{i,t}+\gamma_t(\gamma_tb_{i,t+1}-\beta_t\tilde{\lambda}_{i,t})\right]_+-\tilde{\lambda}_{i,t}\right\|\nonumber\\
&\leq \gamma_t\left\|\gamma_tb_{i,t+1}-\beta_t\tilde{\lambda}_{i,t}\right\|\nonumber\\
&\leq \gamma_t^2\|\nabla g_{i,t}(x_{i,t})(\tilde{x}_{i,t+1}-x_{i,t})+g_{i,t}(x_{i,t})\|+\gamma_t\beta_t\|\tilde{\lambda}_{i,t}\|\nonumber\\
&\leq F\gamma_t^2,\label{equ39}
\end{align} where the first inequality is derived based on $\|[a]_+-[b]_+\|\leq\|a-b\|$ for any two vectors $a,b$ with the same dimension, and the third inequality is obtained based on (\ref{equ11}), (\ref{equ12}) and (\ref{equ28}). Then,
\begin{align}
&\|\lambda_{t+1}-{\bf1}_N\otimes\overline{\lambda}_{t+1}\|\leq\sigma_m\|\lambda_{t}-{\bf1}_N\otimes\overline{\lambda}_t\|+\sqrt{N}F\gamma_t^2.
\end{align}
Combined with $\lambda_{1}-{\bf1}_N\otimes\overline{\lambda}_{1}={\bf0}_{Nm}$, (\ref{equ29}) is thus proved. \hfill$\blacksquare$
\subsection{Proof of Lemma \ref{lemma3}}\label{B}
For any $i\in[N]$, based on the optimality of $\tilde{x}_{i,t+1}$ in (\ref{equ23}), one can obtain that for any $x\in\Omega_i$,
\begin{align}
&\left\langle\tilde{x}_{i,t+1}-x,\alpha_t(\nabla_iJ_{i,t}({\bf x}_{i,t})+\gamma_t(\nabla g_{i,t}(x_{i,t}))^{\top}\tilde{\lambda}_{i,t})\right\rangle\nonumber\\
&+\left\langle\tilde{x}_{i,t+1}-x,\nabla\phi_i(\tilde{x}_{i,t+1})-\nabla\phi_i({x}_{i,t})\right\rangle\leq 0,\label{equ41}
\end{align} where $\frac{\partial D_{\phi_i}(x,x_{i,t})}{\partial x}=\nabla\phi_i({x})-\nabla\phi_i({x}_{i,t})$ for $x\in\Omega_i$ is used. Taking $x=x_{i,t}^*$, it can be derived from (\ref{equ41}) that
\begin{align}
&\alpha_t\left\langle\tilde{x}_{i,t+1}-x_{i,t}^*,\nabla_iJ_{i,t}({\bf x}_{i,t})+\gamma_t(\nabla g_{i,t}(x_{i,t}))^{\top}\tilde{\lambda}_{i,t}\right\rangle\nonumber\\
&\leq \left\langle x_{i,t}^*-\tilde{x}_{i,t+1},\nabla\phi_i(\tilde{x}_{i,t+1})-\nabla\phi_i({x}_{i,t})\right\rangle\nonumber\\
&=D_{\phi_i}(x_{i,t}^*,x_{i,t})-D_{\phi_i}(x_{i,t}^*,\tilde{x}_{i,t+1})-D_{\phi_i}(\tilde{x}_{i,t+1},x_{i,t})\nonumber\\
&\leq D_{\phi_i}(x_{i,t}^*,x_{i,t})-D_{\phi_i}(x_{i,t}^*,\tilde{x}_{i,t+1})\nonumber\\
&~~~-\frac{\mu_0}{2}\|\tilde{x}_{i,t+1}-x_{i,t}\|^2,\label{equ42}
\end{align} where the equality is derived based on (\ref{equ3}) and the last inequality is based on (\ref{equ2}). Rearranging (\ref{equ42}) yields
\begin{align}
&D_{\phi_i}(x_{i,t}^*,\tilde{x}_{i,t+1})\nonumber\\
&\leq\alpha_t\left\langle x_{i,t}^*-\tilde{x}_{i,t+1},\nabla_iJ_{i,t}({\bf x}_{i,t})\right\rangle\nonumber\\
&~~~+\alpha_t\gamma_t\left\langle x_{i,t}^*-\tilde{x}_{i,t+1},(\nabla g_{i,t}(x_{i,t}))^{\top}\tilde{\lambda}_{i,t}\right\rangle\nonumber\\
&~~~+D_{\phi_i}(x_{i,t}^*,x_{i,t})-\frac{\mu_0}{2}\|\tilde{x}_{i,t+1}-x_{i,t}\|^2.\label{equ43}
\end{align} As to the first term of the right hand side of (\ref{equ43}), we have
\begin{align}
&\alpha_t\left\langle x_{i,t}^*-\tilde{x}_{i,t+1},\nabla_iJ_{i,t}({\bf x}_{i,t})\right\rangle\nonumber\\
&=\alpha_t\left\langle x_{i,t}^*-x_{i,t},\nabla_iJ_{i,t}({\bf x}_{i,t})\right\rangle\nonumber\\
&~~~+\alpha_t\left\langle x_{i,t}-\tilde{x}_{i,t+1},\nabla_iJ_{i,t}({\bf x}_{i,t})\right\rangle\nonumber\\
&=\alpha_t\left\langle x_{i,t}^*-x_{i,t},\nabla_iJ_{i,t}(x_{i,t},x_{-i,t})-\nabla_iJ_{i,t}(x_{i,t}^*,x_{-i,t}^*)\right\rangle\nonumber\\
&~~~+\alpha_t\left\langle x_{i,t}^*-x_{i,t},\nabla_iJ_{i,t}({\bf x}_{i,t})-\nabla_iJ_{i,t}(x_{i,t},x_{-i,t})\right\rangle\nonumber\\
&~~~+\alpha_t\left\langle x_{i,t}^*-x_{i,t},\nabla_iJ_{i,t}(x_{i,t}^*,x_{-i,t}^*)+\gamma_t\nabla g_{i,t}(x_{i,t}^*)\lambda_t^*\right\rangle\nonumber\\
&~~~-\alpha_t\gamma_t\left\langle x_{i,t}^*-x_{i,t},\nabla g_{i,t}(x_{i,t}^*)\lambda_t^*\right\rangle\nonumber\\
&~~~+\alpha_t\left\langle x_{i,t}-\tilde{x}_{i,t+1},\nabla_iJ_{i,t}({\bf x}_{i,t})\right\rangle\nonumber\\
&\leq\alpha_t\left\langle x_{i,t}^*-x_{i,t},\nabla_iJ_{i,t}(x_{i,t},x_{-i,t})-\nabla_iJ_{i,t}(x_{i,t}^*,x_{-i,t}^*)\right\rangle\nonumber\\
&~~~+\alpha_t\|x_{i,t}^*-x_{i,t}\|\|\nabla_iJ_{i,t}({\bf x}_{i,t})-\nabla_iJ_{i,t}(x_{i,t},x_{-i,t})\|\nonumber\\
&~~~+\alpha_t\gamma_t\|x_{i,t}^*-x_{i,t}\|\|\nabla g_{i,t}(x_{i,t}^*)\|\|\lambda_t^*\|\nonumber\\
&~~~+\alpha_t\|x_{i,t}-\tilde{x}_{i,t+1}\|\|\nabla_iJ_{i,t}({\bf x}_{i,t})\|\nonumber\\
&\leq\alpha_t\left\langle x_{i,t}^*-x_{i,t},\nabla_iJ_{i,t}(x_{i,t},x_{-i,t})-\nabla_iJ_{i,t}(x_{i,t}^*,x_{-i,t}^*)\right\rangle\nonumber\\
&~~~+2HL\alpha_t\|{\bf x}_{i,t}-x_{t}\|+2LM\Lambda\alpha_t\gamma_t\nonumber\\
&~~~+\frac{M^2\alpha_t^2}{\mu_0}+\frac{\mu_0}{4}\|x_{i,t}-\tilde{x}_{i,t+1}\|^2,\label{equ70}
\end{align} where the first inequality is obtained by Cauchy-Schwarz inequality and the optimality of $x_{i,t}^*$ in (\ref{equ14}), and the second inequality is by Assumption \ref{assump5}, (\ref{equ12}), (\ref{e20}) and Jensen's inequality.

For the second term of the right hand side of (\ref{equ43}), one has
\begin{align}
&\alpha_t\gamma_t\left\langle x_{i,t}^*-\tilde{x}_{i,t+1},(\nabla g_{i,t}(x_{i,t}))^{\top}\tilde{\lambda}_{i,t}\right\rangle\nonumber\\
&=\alpha_t\gamma_t\tilde{\lambda}_{i,t}^{\top}\nabla g_{i,t}(x_{i,t})(x_{i,t}^*-x_{i,t})\nonumber\\
&~~~+\alpha_t\gamma_t\tilde{\lambda}_{i,t}^{\top}\nabla g_{i,t}(x_{i,t})(x_{i,t}-\tilde{x}_{i,t+1})\nonumber\\
&\leq \alpha_t\gamma_t\tilde{\lambda}_{i,t}^{\top}(g_{i,t}(x_{i,t}^*)-g_{i,t}(x_{i,t}))\nonumber\\
&~~~+\alpha_t\gamma_t\tilde{\lambda}_{i,t}^{\top}\nabla g_{i,t}(x_{i,t})(x_{i,t}-\tilde{x}_{i,t+1})\nonumber\\
&= \alpha_t\gamma_t(\tilde{\lambda}_{i,t}-\overline{\lambda}_t)^{\top}(g_{i,t}(x_{i,t}^*)-g_{i,t}(x_{i,t}))\nonumber\\
&~~~+\alpha_t\gamma_t\overline{\lambda}_t^{\top}(g_{i,t}(x_{i,t}^*)-g_{i,t}(x_{i,t}))\nonumber\\
&~~~+\alpha_t\gamma_t\tilde{\lambda}_{i,t}^{\top}\nabla g_{i,t}(x_{i,t})(x_{i,t}-\tilde{x}_{i,t+1})\nonumber\\
&\leq2L\alpha_t\gamma_t\|\tilde{\lambda}_{i,t}-\overline{\lambda}_{t}\|+\alpha_t\gamma_t\overline{\lambda}_{t}^{\top}(g_{i,t}(x_{i,t}^*)-g_{i,t}(x_{i,t}))\nonumber\\
&~~~+\alpha_t\gamma_t\tilde{\lambda}_{i,t}^{\top}\nabla g_{i,t}(x_{i,t})(x_{i,t}-\tilde{x}_{i,t+1})\nonumber\\
&\leq2\sqrt{N}FL\alpha_t\gamma_t\sum\limits_{s=0}^{t-1}\sigma_m^{s}\gamma_{t-1-s}^2+\alpha_t\gamma_t\overline{\lambda}_{t}^{\top}(g_{i,t}(x_{i,t}^*)-g_{i,t}(x_{i,t}))\nonumber\\
&~~~+\alpha_t\gamma_t\tilde{\lambda}_{i,t}^{\top}\nabla g_{i,t}(x_{i,t})(x_{i,t}-\tilde{x}_{i,t+1}),\label{equ71}
\end{align} where the first inequality is derived relaying on the convexity of $g_{ij,t}$ and $\tilde{\lambda}_i\geq0$, the second inequality is obtained by (\ref{equ11}), and the last inequality applies (\ref{equ29}) in Lemma \ref{lemma2}.

By Assumption \ref{assump3} and (\ref{equ24}), one has
\begin{align}
&D_{\phi_i}(x_{i,t}^*,x_{i,t+1})\leq(1-\alpha_t)D_{\phi_i}(x_{i,t}^*,x_{i,t})+\alpha_tD_{\phi_i}(x_{i,t}^*,\tilde{x}_{i,t+1}).\label{equ72}
\end{align} By Assumption \ref{assump2}, it can be obtained that
\begin{align}
&D_{\phi_i}(x_{i,t+1}^*,x_{i,t+1})\leq D_{\phi_i}(x_{i,t}^*,x_{i,t+1})+K\|x_{i,t+1}^*-x_{i,t}^*\|.\label{equ47}
\end{align}

Combining with (\ref{equ43})--(\ref{equ47}) and summing over $i\in[N]$ yields
\begin{align}
&\sum\limits_{i=1}^ND_{\phi_i}(x_{i,t+1}^*,x_{i,t+1})\nonumber\\
&\leq\sum\limits_{i=1}^ND_{\phi_i}(x_{i,t}^*,x_{i,t})+K\sum_{i=1}^N\|x_{i,t+1}^*-x_{i,t}^*\|\nonumber\\
&~~~+2HL\alpha_t^2\sum\limits_{i=1}^N\|{\bf x}_{i,t}-x_{t}\|+2NLM\Lambda\alpha_t^2\gamma_t\nonumber\\
&~~~+\frac{NM^2}{\mu_0}\alpha_t^3+\alpha_t^2\left\langle x_{t}^*-x_{t},F_t(x_t)-F_t(x_t^*)\right\rangle\nonumber\\
&~~~+\alpha_t^2\gamma_t\overline{\lambda}_{t}^{\top}(g_t(x_{t}^*)-g_t(x_{t}))+\alpha_t^2(2c_1(t)-c_2(t)).\label{equ48}
\end{align}

On the other hand, by $\Xi_{t+1}$ in Lemma \ref{lemma2}, one obtains
\begin{align}
-\frac{\Xi_{t+1}}{2\gamma_t}&=-\frac{1}{2}\sum\limits_{i=1}^N\left[\frac{1}{\gamma_t}\|\lambda_{i,t+1}-\lambda\|^2-\frac{1}{\gamma_{t-1}}\|\lambda_{i,t}-\lambda\|^2\right]\nonumber\\
&~~~+\frac{1}{2}\left(\frac{1}{\gamma_t}-\frac{1}{\gamma_{t-1}}-\beta_t\right)\sum_{i=1}^N\|\lambda_{i,t}-\lambda\|^2. \label{equ49}
\end{align}
Summing over $t\in[T]$ gives that
\begin{align}
&-\sum_{t=1}^T\frac{\Xi_{t+1}}{2\gamma_t}\nonumber\\
&=-\frac{1}{2}\sum\limits_{i=1}^N\left[\frac{1}{\gamma_T}\|\lambda_{i,T+1}-\lambda\|^2-\frac{1}{\gamma_{0}}\|\lambda_{i,1}-\lambda\|^2\right]\nonumber\\
&~~~+\frac{1}{2}\sum\limits_{t=1}^T\left(\frac{1}{\gamma_t}-\frac{1}{\gamma_{t-1}}-\beta_t\right)\sum_{i=1}^N\|\lambda_{i,t}-\lambda\|^2\nonumber\\
&\leq\frac{N}{2}\|\lambda\|^2+\frac{1}{2}\sum\limits_{t=1}^T\left(\frac{1}{\gamma_t}-\frac{1}{\gamma_{t-1}}-\beta_t\right)\sum_{i=1}^N\|\lambda_{i,t}-\lambda\|^2,\label{equ50}
\end{align} where the inequality is derived based on $\lambda_{i,1}={\bf0}_m$.

Based on Assumption \ref{assump6}, summing over $t\in[T]$ on both sides of (\ref{equ48}), combining with (\ref{equ50}) and (\ref{equ30}) in Lemma \ref{lemma2}, we have
\begin{align}
&\mu\sum_{t=1}^T\|x_t-x_t^*\|^2\nonumber\\
&\leq\sum_{t=1}^T\left\langle x_{t}^*-x_{t},F_t(x_t^*)-F_t(x_t)\right\rangle\nonumber\\
&\leq\sum\limits_{t=1}^T\frac{1}{\alpha_t^2}\sum_{i=1}^N\left[D_{\phi_i}(x_{i,t}^*,x_{i,t})-D_{\phi_i}(x_{i,t+1}^*,x_{i,t+1})\right]\nonumber\\
&~~~+\frac{\sqrt{N}K}{\alpha_T^2}\Phi_T^*+3\sum_{t=1}^Tc_1(t)+2HL\sum_{t=1}^T\sum\limits_{i=1}^N\|{\bf x}_{i,t}-x_{t}\|\nonumber\\
&~~~+2NLM\Lambda\sum_{t=1}^T\gamma_t+\frac{NF^2}{2}\sum_{t=1}^T\gamma_t^3+\frac{NM^2}{\mu_0}\sum_{t=1}^T\alpha_t\nonumber\\
&~~~+\sum_{t=1}^T\gamma_t\overline{\lambda}_t^{\top}g_t(x_t^*)-\lambda^{\top}\sum_{t=1}^T\gamma_tg_t(x_t)+\frac{N}{2}\|\lambda\|^2\nonumber\\
&~~~+\frac{1}{2}\sum\limits_{t=1}^T\left(\frac{1}{\gamma_t}-\frac{1}{\gamma_{t-1}}-\beta_t\right)\sum_{i=1}^N\|\lambda_{i,t}-\lambda\|^2\nonumber\\
&~~~+N\sum_{t=1}^T\left(\frac{M^2}{\mu_0}\alpha_t\gamma_t^2+\frac{1}{2}\beta_t\right)\|\lambda\|^2. \label{equ51}
\end{align}

Because
\begin{align}
\sum_{t=1}^T\sum\limits_{s=0}^{t-1}\sigma^s\gamma_{t-s-1}&=\sum_{t=1}^T\gamma_{t-1}\sum\limits_{s=0}^{T-t}\sigma^s\nonumber\\
&\leq\frac{1}{1-\sigma}\sum_{t=1}^T\gamma_{t-1},\label{equ52}
\end{align}
we get
\begin{align}
\sum_{t=1}^T\gamma_t\sum_{s=0}^{t-1}\sigma_m^{s}\gamma_{t-1-s}^2&\leq\sum_{t=1}^T\sum_{s=0}^{t-1}\sigma_m^{s}\gamma_{t-1-s}^3\nonumber\\
&\leq\frac{1}{1-\sigma_m}\sum_{t=1}^T\gamma_{t-1}^3,\label{equ53}\\
\sum_{t=1}^T\sum_{i=1}^N\|{\bf x}_{i,t}-x_t\|&\leq\sqrt{N}\sum_{t=1}^T\sqrt{\sum_{i=1}^N\|{\bf x}_{i,t}-x_t\|^2}\nonumber\\
&\leq\sqrt{N}\sum_{t=1}^T\sum_{i=1}^N\|e_{i,t}\|\nonumber\\
&\leq2N^2L\sum_{t=1}^T\sum\limits_{s=0}^{t-1}\sigma^s\alpha_{t-s-1}\nonumber\\
&\leq\frac{2N^2L}{1-\sigma}\sum_{t=1}^T\alpha_{t-1},
\label{equ54}
\end{align} where the third inequality is derived based on Lemma \ref{lemma1}. Substituting (\ref{equ52})--(\ref{equ54}) into (\ref{equ51}) and in view of $g_t(x_t^*)\leq{\bf0}_m$ together with $\gamma_{t}^3\leq\gamma_t$, one can conclude that (\ref{equ34}) holds. \hfill$\blacksquare$

\subsection{Proof of Lemma \ref{lemma4}}\label{C}
For the first term on the right hand side of (\ref{equ34}), it can be derived that
{\small\begin{align}
&\sum_{t=1}^T\frac{1}{\alpha_t^2}\sum_{i=1}^N\left[D_{\phi_i}(x_{i,t}^*,x_{i,t})-D_{\phi_i}(x_{i,t+1}^*,x_{i,t+1})\right]\nonumber\\
&\leq\sum_{t=1}^T\sum_{i=1}^N\left[\frac{1}{\alpha_t^2}D_{\phi_i}(x_{i,t}^*,x_{i,t})-\frac{1}{\alpha_{t+1}^2}D_{\phi_i}(x_{i,t+1}^*,x_{i,t+1})\right]\nonumber\\
&~~~+\sum_{t=1}^T\sum_{i=1}^N\left[\frac{1}{\alpha_{t+1}^2}-\frac{1}{\alpha_t^2}\right]D_{\phi_i}(x_{i,t+1}^*,x_{i,t+1})\nonumber\\
&\leq\frac{1}{\alpha_1^2}\sum_{i=1}^ND_{\phi_i}(x_{i,1}^*,x_{i,1})+\sum_{t=1}^T\sum_{i=1}^N\left[\frac{1}{\alpha_{t+1}^2}-\frac{1}{\alpha_t^2}\right]2LK\nonumber\\
&\leq\frac{2NLK}{\alpha_{T+1}^2},\label{equ56}
\end{align}}where Assumptions \ref{assump2} and \ref{assump4} have been used to get the second inequality. Then, setting $\lambda={\bf 0}_{m}$, (\ref{equ34}) becomes
\begin{align}
&\mu\sum_{t=1}^T\|x_{t}-x_{t}^*\|^2\nonumber\\
&\leq\frac{2NLK}{\alpha_{T+1}^2}+\frac{\sqrt{N}K}{\alpha_T^2}\Phi_T^*+C_1\sum_{t=1}^T\gamma_{t-1}+C_2\sum_{t=1}^T\alpha_t\nonumber\\
&~~~+\frac{1}{2}\sum_{t=1}^T\sum\limits_{i=1}^N\left(\frac{1}{\gamma_t}-\frac{1}{\gamma_{t-1}}-\beta_t\right)\|\lambda_{i,t}\|^2.\label{equ55}
\end{align}
Together with
\begin{align}
Reg_i(T)&=\sum_{t=1}^T(J_{i,t}(x_{i,t},x_{-i,t}^*)-J_{i,t}(x_{i,t}^*,x_{-i,t}^*))\nonumber\\
&\leq M\sum_{t=1}^T\|x_{i,t}-x_{i,t}^*\|\nonumber\\
&\leq M\sqrt{T\sum_{t=1}^T\|x_{t}-x_{t}^*\|^2},
\end{align} one can conclude that (\ref{e35}) holds.

Let $\lambda=\lambda_c=\frac{2\left[\sum_{t=1}^T\gamma_tg_t(x_t)\right]_+}{B_3(T)}$, then (\ref{equ34}) and (\ref{equ56}) imply that
\begin{align}
&\mu\sum_{t=1}^T\|x_{t}-x_{t}^*\|^2\nonumber\\
&\leq \mu B_1(T)+\frac{\sqrt{N}K}{\alpha_T^2}\Phi_T^*-\frac{1}{B_3(T)}\left\|\left[\sum\limits_{t=1}^T\gamma_tg_t(x_t)\right]_+\right\|^2\nonumber\\
&~~~+\frac{1}{2}\sum\limits_{t=1}^T\left(\frac{1}{\gamma_t}-\frac{1}{\gamma_{t-1}}-\beta_t\right)\sum_{i=1}^N\|\lambda_{i,t}-\lambda_c\|^2,\label{e60}
\end{align} from which, it can be derived that
\begin{align}
&\left\|\left[\sum\limits_{t=1}^T\gamma_tg_t(x_t)\right]_+\right\|^2\nonumber\\
&\leq-\mu B_3(T)\sum_{t=1}^T\|x_{t}-x_{t}^*\|^2+\mu B_1(T)B_3(T)+\frac{\sqrt{N}K}{\alpha_T^2}B_3(T)\Phi_T^*\nonumber\\
&~~~+\frac{B_3(T)}{2}\sum\limits_{t=1}^T\left(\frac{1}{\gamma_t}-\frac{1}{\gamma_{t-1}}-\beta_t\right)\sum_{i=1}^N\|\lambda_{i,t}-\lambda_c\|^2.\label{e61}
\end{align} Based on $\|x_t-x_t^*\|^2=\sum_{i=1}^N\|x_{i,t}-x_{i,t}^*\|^2\leq4NL^2$, Lemma \ref{lemma4} is thus proved. \hfill$\blacksquare$

\end{appendix}
\bibliographystyle{IEEEtran}
\bibliography{MinMeng}

\end{document}